\documentclass{amsart}
\usepackage{amssymb, mathrsfs}
\usepackage[all]{xy}

\newcommand{\Bc}[9]{\bibitem{#1} {#2}, \emph{#3}, in: \textbf{#4} (#5), #6 #7, #8--#9.}
\newcommand{\nInc}{\compactN^{\uparrow n}}
\newcommand{\Increasing}{{\compactN^{\uparrow\N}}}
\newcommand{\Inc}{\Increasing}
\newcommand{\powInc}[1]{{\Increasing^{#1}}}
\newcommand{\gpbl}{{\mbox{\textit{\tiny gp}}}}
\newcommand{\SMZ}{\mathrm{SMZ}}
\newcommand{\sub}{\op{sub}}

\long\def\forget#1\forgotten{}

\newcommand{\cU}{\mathcal{U}}
\newcommand{\cV}{\mathcal{V}}
\newcommand{\cP}{\mathcal{P}}

\newcommand{\Dfin}{\mathfrak{D}_\mathrm{fin}}

\newcommand{\G}{\mathcal{G}}
\newcommand{\e}{\max\{\s,\b\}}
\newcommand{\sr}[2]{{\txt{{${#1}$}\\{\tiny ${#2}$}}}}
\newcommand{\scrA}{\mathscr{A}}
\newcommand{\scrB}{\mathscr{B}}

\newcommand{\seq}[1]{\{#1\}_{n\in\N}}

\newcommand{\op}{\operatorname}

\newcommand{\Cantor}{{{}^\N\{0,1\}}}
\newcommand{\Cal}{\mathcal}

\newcommand{\B}{{\Cal B}}

\newcommand{\BG}{\B_\Gamma}
\newcommand{\BT}{\B_\mathrm{T}}
\newcommand{\BO}{\B_\Omega}

\newcommand{\Tau}{\mathrm{T}}

\newcommand{\cF}{{\Cal F}}
\newcommand{\M}{\mathcal{M}}
\newcommand{\N}{\mathbb{N}}
\newcommand{\compactN}{\cl{\mathbb{N}}}

\newcommand{\NN}{{{}^\N\N}}
\newcommand{\NcompactN}{{{}^\N\compactN}}

\newcommand{\roth}{{[\N]^{\aleph_0}}} 
\newcommand{\fin}{{[\N]^{<\aleph_0}}} 
\newcommand{\cO}{\mathcal{O}}
\newcommand{\Q}{\mathbb{Q}}
\newcommand{\R}{\mathbb{R}}

\newcommand{\U}{{\Cal U}}
\newcommand{\Union}{\bigcup}

\renewcommand{\b}{\mathfrak{b}}
\renewcommand{\c}{\mathfrak{c}}
\renewcommand{\d}{\mathfrak{d}}
\newcommand{\fu}{\mathfrak{u}}

\newcommand{\g}{\mathfrak{g}}
\newcommand{\oo}{\infty}
\newcommand{\p}{\mathfrak{p}}
\newcommand{\s}{\mathfrak{s}}
\newcommand{\w}{\omega}
\newcommand{\x}{\times}

\newcommand{\nin}{\not\in}

\newcommand{\sbst}{\subseteq}
\newcommand{\spst}{\supseteq}
\newcommand{\sm}{\setminus}

\newcommand{\as}{\sbst^*}

\newcommand{\J}{\mathcal{J}}

\newcommand{\cov}{{\sf cov}}
\newcommand{\add}{{\sf add}}

\newcommand{\non}{{\sf non}}

\renewcommand{\t}{\mathfrak{t}}

\newtheorem{thm}{Theorem}
\newtheorem{prop}[thm]{Proposition}
\newtheorem{prob}[thm]{Problem}

\newtheorem{lem}[thm]{Lemma}
\newtheorem{cor}[thm]{Corollary}

\theoremstyle{definition}

\theoremstyle{remark}
\newtheorem{rem}[thm]{Remark}

\newcommand{\be}{\begin{enumerate}}
\newcommand{\ee}{\end{enumerate}}
\newcommand{\bi}{\begin{itemize}}
\newcommand{\ei}{\end{itemize}}

\newcommand{\sone}{{\sf S}_1}    \newcommand{\sfin}{{\sf S}_{fin}}
\newcommand{\ufin}{{\sf U}_{fin}}

\newcommand{\cl}{\overline}

\newcommand{\rest}{{\mathord{\restriction}}}
\newcommand{\lft}[2]{\mathopen\ifcase#1{}\oo\or
                        \big#2\or\Big#2\else\oo\fi}
\newcommand{\rgt}[2]{\mathclose\ifcase#1{}\oo\or
                        \big#2\or\Big#2\else\oo\fi}

\title[Hereditary topological diagonalizations]{Hereditary topological diagonalizations
and the Menger-Hurewicz Conjectures}
\author{Tomek Bartoszy\'nski}
\thanks{The first author was  partially supported by
NSF grant DMS 0200671.}
\address{Department of Mathematics and Computer Science,
Boise State University, Boise, Idaho 83725 U.S.A.
}
\email{tomek@math.boisestate.edu}
\urladdr{http://math.boisestate.edu/\~{}tomek}
\author{Boaz Tsaban}
\thanks{This paper constitutes a part of the second author's doctoral dissertation at
Bar-Ilan University.}
\address{Department of Applied Mathematics and Computer Science,
Weizmann Institute of Science, Rehovot 76100, Israel}
\email{boaz.tsaban@weizmann.ac.il}
\urladdr{http://www.cs.biu.ac.il/\~{}tsaban}

\begin{document}
\begin{abstract}
We consider the question, which of the major classes defined by
topological diagonalizations of open or Borel covers is hereditary.
Many of the classes in the open case are not hereditary already in ZFC,
and none of them is provably hereditary. This is contrasted with the
Borel case, where some of the classes are provably hereditary.
Two of the examples are counter-examples of sizes $\d$ and $\b$,
respectively, to the Menger and Hurewicz Conjectures, and one of
them answers a question of Steprans on perfectly meager sets.

\medskip

\emph{Remark.} The main results of this paper are improved in the paper
\emph{Menger's and Hurewicz's Problems: Solutions from ``The Book'' and refinements},\\
\texttt{http://arxiv.org/abs/0909.5645}
\end{abstract}


\maketitle

\section{Introduction}

\subsection{Selection principles}
Let $\scrA$ and $\scrB$ be collections of covers of a topological
space $X$. The following selection hypotheses have a long history
for the case where the collections $\scrA$ and $\scrB$ are
topologically significant.
\begin{itemize}
\item[$\sone(\scrA,\scrB)$:]{
For each sequence $\seq{\cU_n}$ of members of $\scrA$,
there exist members $U_n\in\cU_n$, $n\in\N$, such that $\seq{U_n}\in\scrB$.
}
\item[$\sfin(\scrA,\scrB)$:]{
For each sequence $\seq{\cU_n}$
of members of $\scrA$, there exist finite (possibly empty)
subsets $\cF_n\sbst\cU_n$, $n\in\N$, such that $\Union_{n\in\N}\cF_n\in\scrB$.
}
\item[$\ufin(\scrA,\scrB)$:]{
For each sequence $\seq{\cU_n}$ of members of $\scrA$
which do not contain a finite subcover,
there exist finite (possibly empty) subsets $\cF_n\sbst\cU_n$, $n\in\N$,
such that $\seq{\cup\cF_n}\in\scrB$.
}
\end{itemize}

\subsection{Special covers}
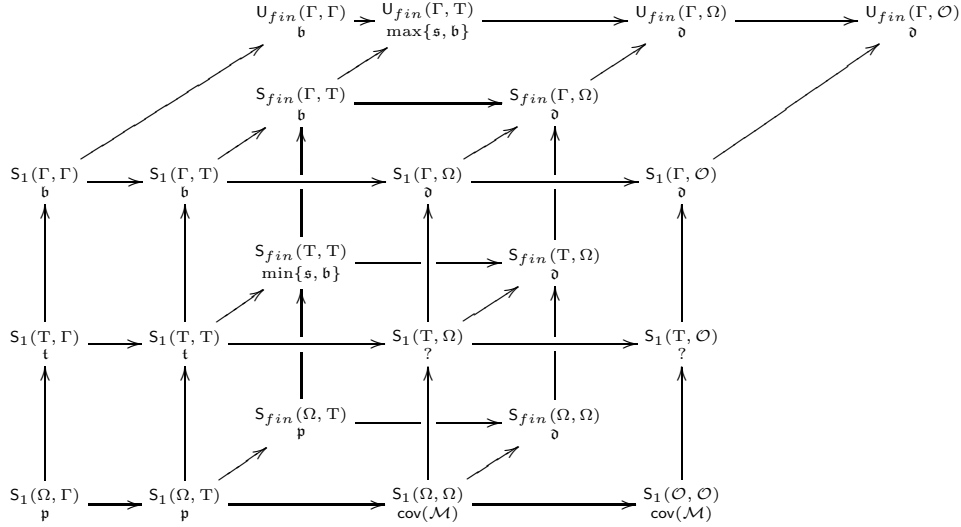
\begin{figure}[!ht]
{\tiny
$$\xymatrix@C=7pt@R=12pt{
&
&
& \sr{\ufin(\Gamma,\Gamma)}{\b}\ar[r]
& \sr{\ufin(\Gamma,\Tau)}{\e}\ar[rr]
&
& \sr{\ufin(\Gamma,\Omega)}{\d}\ar[rrrr]
&
&
&
& \sr{\ufin(\Gamma,\cO)}{\d}
\\
&
&
& \sr{\sfin(\Gamma,\Tau)}{\b}\ar[rr]\ar[ur]
&
& \sr{\sfin(\Gamma,\Omega)}{\d}\ar[ur]
\\
\sr{\sone(\Gamma,\Gamma)}{\b}\ar[uurrr]\ar[rr]
&
& \sr{\sone(\Gamma,\Tau)}{\b}\ar[ur]\ar[rr]
&
& \sr{\sone(\Gamma,\Omega)}{\d}\ar[ur]\ar[rr]
&
& \sr{\sone(\Gamma,\cO)}{\d}\ar[uurrrr]
\\
&
&
& \sr{\sfin(\Tau,\Tau)}{\min\{\s,\b\}}\ar'[r][rr]\ar'[u][uu]
&
& \sr{\sfin(\Tau,\Omega)}{\d}\ar'[u][uu]
\\
\sr{\sone(\Tau,\Gamma)}{\t}\ar[rr]\ar[uu]
&
& \sr{\sone(\Tau,\Tau)}{\t}\ar[uu]\ar[ur]\ar[rr]
&
& \sr{\sone(\Tau,\Omega)}{\mathbf{?}}\ar[uu]\ar[ur]\ar[rr]
&
& \sr{\sone(\Tau,\cO)}{\mathbf{?}}\ar[uu]
\\
&
&
& \sr{\sfin(\Omega,\Tau)}{\p}\ar'[u][uu]\ar'[r][rr]
&
& \sr{\sfin(\Omega,\Omega)}{\d}\ar'[u][uu]
\\
\sr{\sone(\Omega,\Gamma)}{\p}\ar[uu]\ar[rr]
&
& \sr{\sone(\Omega,\Tau)}{\p}\ar[uu]\ar[ur]\ar[rr]
&
& \sr{\sone(\Omega,\Omega)}{\cov(\M)}\ar[uu]\ar[ur]\ar[rr]
&
& \sr{\sone(\cO,\cO)}{\cov(\M)}\ar[uu]
}$$
}
\caption{The surviving classes for the open case}\label{survopen}
\end{figure}

Let $X$ be a
set of reals. In the following definitions, we always require that
$X$ is not contained in any member of the cover. An
\emph{$\w$-cover} $\U$ of $X$ is a cover of $X$ such that each
finite subset of $X$ is contained in some member of $\U$. $\U$ is
a \emph{$\tau$-cover} if it is large (that is, each member of $X$
is contained in infinitely many members of the cover), and for
each $x,y\in X$, (at least) one of the sets $\{U\in\U : x\in U,
y\nin U\}$ and $\{U\in\U : y\in U, x\nin U\}$ is finite. $\U$ is a
$\gamma$-cover if it is infinite, and each element of $X$ belongs
to all but finitely many members of the cover. Let $\cO$, $\Omega$,
$\Tau$, and $\Gamma$ denote the collections of countable open
covers, $\w$-covers, $\tau$-covers, and $\gamma$-covers of $X$,
respectively, and let $\B,\BO,\BT,\BG$ be the corresponding
countable \emph{Borel} covers. The diagonalization properties of
these types of covers were extensively studied in, e.g.,
\cite{comb1, coc2, CBC, tautau}. Many of
these properties turn out equivalent \cite{coc2, tautau}; the
classes which survived thus far appear in Figure~\ref{survopen}
(for the open case). Some of the classes which are distinct in the
open case coincide in the Borel case \cite{CBC}.

In the diagram, each property $P$ appears together with its
\emph{critical cardinality} $\non(P)$,
that is, the minimal size of a set of reals which does not
satisfy that property.
(See \cite{vD, HBK} for the definitions of these and other constants
mentioned in the paper.)
The arrows in this diagram denote inclusion.

\subsection{Hereditary properties}
If $\scrA$ and $\scrB$ are collections of open covers and
$\Pi\in\{\sone,\sfin,\ufin\}$, then $\Pi(\scrA,\scrB)$ is closed
under taking closed subsets \cite{coc2}. Similarly, if $\scrA$ and
$\scrB$ are collections of Borel covers and
$\Pi\in\{\sone,\sfin,\ufin\}$, then $\Pi(\scrA,\scrB)$ is closed
under taking Borel subsets \cite{CBC}. We say that a class (or a
property) is \emph{hereditary} if it is closed under taking
\emph{arbitrary} subsets. A natural question which rises is: Which
of these classes is provably hereditary? We show that some of the
classes in the Borel case, but none of the classes in the open
case, are provably hereditary.

For some of the classes in the open case,
no additional axioms beyond ZFC are required in order to
disprove their being hereditary.
We describe two nontrivial examples which are in particular counter-examples to
conjectures of Menger and Hurewicz.

\section{The Borel case}

\begin{prop}\label{BBhered}
$\sone(\B,\B)$ is hereditary.
\end{prop}
\begin{proof}
Assume that $X$ satisfies $\sone(\B,\B)$ and that $Y$ is a subset of $X$.
Assume that $\seq{\cU_n}$ is a sequence of countable Borel covers of $Y$.

For each $n$, define $\cV_n = \cU_n\cup\{X\sm\cup\cU_n\}$.
Then each $\cV_n$ is a countable Borel cover of $X$, and we can choose
for each $n$ an element $V_n\in\cV_n$ such that $\seq{V_n}$ is a
cover of $X$.
Define $U_n = V_n$ if $V_n\neq X\sm\cup\cU_n$, otherwise let $U_n$ be an arbitrary
element of $\cU_n$. Then $\seq{U_n}$ is a cover of $Y$ and for each $n$,
$U_n\in\cU_n$.
\end{proof}

All classes in the Borel case of the top plane of Figure~\ref{survopen}
are hereditary.
\begin{prop}\label{referee}
For each $\Pi\in\{\sone,\sfin,\ufin\}$ and each
$\scrB\in\{\B,\BO,\BT,\BG\}$, $\Pi(\BG,\scrB)$ is hereditary.
\end{prop}
\begin{proof}
The proof for this is similar to that of Proposition \ref{BBhered}:
If $\cU= \seq{U_n}$ is a countable Borel $\gamma$-cover of a subset $Y$ of $X$,
then
$$B_\cU = \{x\in X : \mbox{for infinitely many $n$, }x\nin U_n\}$$
is a Borel subset of $X$ disjoint from $Y$.
We claim that $\cV=\seq{U_n\cup B_\cU}$ is a (countable Borel) $\gamma$-cover of $X$.
It is easy to see that each $x\in X$ is contained
in all but finitely many members of $\cV$.
It remains to show that the cover is infinite.
As $Y$ is not contained in any element of $\cU$ and $B_\cU$ is disjoint from
$Y$, $X$ is not contained in any member of $\cV$.
Moreover, each finite subset of $X$ is contained in some element
of $\cV$. Thus, $\cV$ is an $\omega$-cover of $X$, and in particular it is infinite.
\end{proof}

\begin{prob}
Which of the remaining classes (which involve Borel covers)
are provably hereditary?
\end{prob}

Recently, Miller answered a part of this problem by showing that
no class between $\sone(\BO,\BG)$ and $\sfin(\Omega,\Tau)$ (inclusive)
is provably hereditary \cite{MilNonGamma}.

\subsection*{Consistency facts}\label{con}
Borel's Conjecture, which was proved consistent by Laver,
implies that the class $\sone(\cO,\cO)$ (and the classes below it) reduces to contain
only the countable sets of reals.
Thus, it is consistent that all classes below $\sone(\cO,\cO)$ are hereditary.

A set of reals $X$ is a \emph{$\sigma$-set} if each $G_\delta$ set in $X$
is also an $F_\sigma$ set in $X$.
By \cite{CBC}, every element of $\sone(\BG,\BG)$
is a $\sigma$-set.
It is consistent that every $\sigma$-set of
real numbers is countable \cite{Mil79Len}.
Consequently, it is consistent that all classes below $\sone(\BG,\BG)$ are hereditary.

It is a major open problem whether it is consistent that every uncountable set of
real numbers can be mapped onto a \emph{dominating} subset of $\NN$ by a Borel function.
Such a consistency result would imply the consistency of all classes
considered in this paper for the Borel case being hereditary.

\section{The open case}
A quasiordering $\le^*$ is defined on the Baire space $\NN$:
$f\le^* g$ if $f(n)\le g(n)$ for all but finitely many $n$.
A subset $Y$ of $\NN$ is \emph{dominating} if for each $g$ in $\NN$ there
exists $f\in Y$ such that $g\le^* f$. It is \emph{bounded} if
there exists $g\in\NN$ such that for each $f\in Y$, $f\le^* g$.
\emph{Cantor's space} $\Cantor$ of infinite
binary sequences is equipped with the product topology.
Identify $\Cantor$ with $P(\N)$ by characteristic functions.
Also, denote by $\roth$ (respectively, $\fin$) the collection of all
infinite (respectively, finite) sets of natural numbers.

Assuming Martin's Axiom (or just $\p=\c$),
there exists $X\sbst P(\N)$ of size $\c$
such that $X$ satisfies $\sone(\Omega,\Gamma)$
but $X\sm\fin$ does not satisfy $\sone(\Omega,\Gamma)$ \cite{GM}.
We will modify the construction of \cite{GM} to get a stronger result.
A space $X$ is a \emph{$\tau$-set} \cite{tau}
if each clopen $\tau$-cover of $X$ contains
a $\gamma$-cover of $X$.

\begin{thm}\label{p=c-X}
Assume that $\p=\c$. Then there exists $X\sbst P(\N)$
such that $X$ satisfies $\sone(\Omega,\Gamma)$
but $X\sm\fin$ does not satisfy $\ufin(\Gamma,\cO)$.
(Moreover, $X$ is a $\tau$-set and $X\sm\fin$ is not a $\tau$-set.)
\end{thm}

\begin{proof}
For $y\sbst\N$, define $y^*=\{x\sbst\N : x\as y\}$.

\begin{lem}[\cite{GM}]\label{GMlemma}
Assume that $\G$ is an open $\w$-cover of $\fin$.
Then for each infinite $x\sbst\N$ there exists an infinite
$y\sbst x$ such that
$\G$ contains a $\gamma$-cover of $y^*$.
\end{lem}

Identify $P(\N)$ with the collection of strictly increasing
functions in $\NN$ by taking increasing enumerations.
When $a\in\roth$, $\vec a$ denotes the increasing enumeration of $a$.

Let $\{\G_\alpha : \alpha<\c\}$ enumerate all countable families of open
sets in $P(\N)$, $\{\vec d_\alpha : \alpha<\c\}$ be a dominating subset of $\NN$,
and $\{a_\alpha : \alpha<\c\}\sbst\roth$ enumerate all infinite coinfinite
subsets of $\N$.
For convenience, for $x,y\in\roth$ we define
$$\sub(x,y)=\begin{cases}
x      & x\as y\\
x\sm y & \mbox{otherwise}
\end{cases}$$
Thus $\sub(x,y)\in\roth$,
$\sub(x,y)\sbst x$,
and either $\sub(x,y)\as y$ or else $\sub(x,y)\as\N\sm y$.

We construct by induction a dominating subset
$\{\vec x_\alpha : \alpha<\c\}$ of $\NN$, such that
$\{x_\alpha : \alpha<\c\}\sbst\roth$ is a (special type of a)
tower, as follows.

For a limit $\alpha$, use $\alpha<\c=\b=\t$ to get a
pseudo intersection $p$ of $\{x_\beta : \beta<\alpha\}$ and
a function $b\in\NN$ which bounds $\{\vec d_\beta : \beta<\alpha\}$.
Choose an infinite $q\sbst p$ such that $\vec b\le^* \vec q$,
e.g., $\vec q(n)=\min\{k\in p : \vec q(n-1),\vec b(n)<k \}$.
Then set $x_\alpha = \sub(q,a_\alpha)$.

The successors $x_{\alpha+1}$ are constructed as follows.
If $\G_\alpha$ is not an $\omega$-cover of
$X_\alpha = \{x_\beta : \beta<\alpha\}\cup\fin$, set
$x_{\alpha+1}=\sub(x_\alpha,a_{\alpha+1})$
(this case is not particularly interesting). Otherwise
do the following:
As $|X_\alpha|<\p$, $\G_\alpha$ contains a $\gamma$-cover of
$X_\alpha$.
By Lemma \ref{GMlemma},
there exists an infinite $p\sbst x_\alpha$ such that
this $\gamma$-cover (which is in particular an $\w$-cover of $\fin$)
contains a $\gamma$-cover $\seq{G_n}$ of $p^*$. Observe that
$\seq{G_n}$ is a $\gamma$-cover of $ \{x_\beta : \beta<\alpha\}\cup p^*$.
Now, as in the first case,
take an infinite $q\sbst p$
such that $\vec q$ bounds $\{\vec d_\beta : \beta<\alpha+1\}$,
and set $x_{\alpha+1} = \sub(q,a_{\alpha+1})$.

Set $X=\{x_\alpha : \alpha<\c\}\cup\fin$.
The properties follow, as in
Theorems 2.14 and 3.7 of \cite{tau}. Briefly:
To prove that $X$ satisfies $\sone(\Omega,\Gamma)$, it suffices to show that each
$\w$-cover of $X$ contains a $\gamma$-cover of $X$ \cite{GN}.
By the construction, if $\G_\alpha$ is an $\omega$-cover of
$X$, then it contains a $\gamma$-cover
of $\{x_\beta : \beta\le \alpha\}\cup x_{\alpha+1}^*\spst X$.

Each $\gamma$-set is a $\tau$-set. However, a tower is never a $\tau$-set
\cite{tau}. $X\sm\fin$ is a tower:
Let $a\in\roth$.
Take an infinite coinfinite $a_\alpha\sbst a$.
Then either $x_\alpha\as a_\alpha$, or else $x_\alpha\as \N\sm a_\alpha$.
Therefore, $a\not\as x_\alpha$.

Finally, by a theorem of Hurewicz, 
a zero-dimensional space $X$ satisfies $\ufin(\Gamma,\cO)$ if, and only if,
for every continuous function $\Psi:X\to\NN$, $\Psi[X]$ is not dominating.
As $X\sm\fin$ is dominating, it does not satisfy $\ufin(\Gamma,\cO)$.
\end{proof}

A construction as in Theorem \ref{p=c-X} cannot give information
in the Borel cases: E.g.,
if $X$ satisfies $\sone(\BO,\BG)$ and $A$ is countable, then
$X\sm A$ is a Borel subset of $X$ and therefore satisfies $\sone(\BO,\BG)$ too.
However, this construction can be strengthened to obtain a
\emph{strong} $\gamma$-set \cite{GM}.
A set $X$ of reals is a \emph{strong $\gamma$-set} if, and only if,
for each sequence $\seq{\U_n}$, where for each $n$ $\U_n$ is an $n$-cover of $X$,
there exist elements $U_n\in\U_n$, $n\in\N$, such that $\seq{U_n}$
is a $\gamma$-cover of $X$ \cite{strongdiags}.

\begin{thm}
Assume Martin's Axiom. Then there exists $X\sbst P(\N)$
such that $X$ is a strong $\gamma$-set
but $X\sm\fin$ does not satisfy $\ufin(\Gamma,\cO)$ and
is not a $\tau$-set.
\end{thm}
\begin{proof}
We carry out the same construction as in Theorem
\ref{p=c-X}, but replace Lemma \ref{GMlemma}
with the following one, which is also due to Galvin and Miller \cite{GM}.
\begin{lem}[MA]
Assume that $X\sbst P(\N)$ is such that $|X|<\c$,
and $x\in\roth$.
Then for each sequence $\seq{\U_n}$ of open $n$-covers of $X\cup\fin$
there exists an infinite subset $y$ of $x$ and a sequence
$\seq{V_n}$ such that for each $n$ $V_n\in\U_n$, and $\seq{V_n}$
is a $\gamma$-cover of $X\cup y^*$.
\end{lem}
This allows us to carry out the construction where we consider
all possible sequences of $n$-covers instead of all possible $\omega$-covers.
\end{proof}

\section{The Menger and Hurewicz conjectures}

Some of the classes can be treated without any special hypotheses.

\begin{lem}
Assume that $\J\sbst P(\R)$ is closed
under taking subsets and continuous images, and
assume that $[0,1]\in\J$. Then $\J=P(\R)$.
\end{lem}
\begin{proof}
The subset $(0,1)$ of $[0,1]$ belongs to $\J$ and can be mapped
continuously onto $\R$. Thus, every subset $X$ of $\R$ is a
continuous image of some subset of $(0,1)$.
\end{proof}

By the definition, every $\sigma$-compact set of reals
satisfies $\ufin(\Gamma,\Gamma)$.
Moreover, every $\sigma$-compact set satisfies $\sfin(\Omega,\Omega)$ \cite{coc2}.
Since all of the properties in Figure~\ref{survopen} are closed under
taking continuous images and $\ufin(\Gamma,\cO)\neq P(\R)$,
we have the following.
\begin{cor}
None of the classes in Figure~\ref{survopen} which contain
$\ufin(\Gamma,\Gamma)$ or $\sfin(\Omega,\Omega)$ is hereditary.
\end{cor}

A natural question is whether all nonhereditary sets in these
classes are $\sigma$-compact.
Related questions were raised by Menger and Hurewicz.
Menger \cite{MENGER} conjectured that each set of reals satisfying
$\ufin(\Gamma,\cO)$ is $\sigma$-compact.
Hurewicz \cite{HURE25} had a weaker conjecture that
$\ufin(\Gamma,\Gamma)$ implies $\sigma$-compactness.
A Sierpinski set is a counter-example to both conjectures \cite{coc2},
but it was only recently \cite{FM, coc2} that these conjecture were
disproved in ZFC.
Both refutations use a dichotomic argument: For an appropriate cardinal
$\kappa\geq\aleph_1$, two independent examples are given
for the case $\kappa=\aleph_1$ and for the case $\kappa>\aleph_1$.
Our goal in the sequel is to explore the properties of two axiomatic-independent
counter-examples to these conjectures.

We adopt the following setting from \cite{BaShCon2000}.
Let $\compactN=\N\cup\{\infty\}$ be the one point compactification of $\N$.
A subset $A\sbst\compactN$ is open if:
$A\sbst\N$, or
$\infty\in A$ and $A$ is cofinite.
Thus, if $A$ is a compact subset of $\compactN$ and $\infty\nin A$, then
$A$ is finite.

Let $\Increasing\sbst\NcompactN$ consist of the nondecreasing functions $f$ such
that for each $n$ with $f(n)<\infty$, $f(n)<f(n+1)$.
$\Increasing$ is a zero-dimensional metrizable compact space
without isolated points, thus by a classical theorem of Brouwer it is
homeomorphic to the Cantor space $\Cantor$.
For each increasing finite sequence $s$ of natural numbers,
let $q_s\in \Increasing$ be defined as
$$q_{s}(k)=
\begin{cases}
s(k) & \text{if } k <|s|\\
\infty & \text{otherwise}
\end{cases}$$
for each $k\in\N$.
The set $Q$ of all these elements $q_s$ is dense in $\Increasing$.

The following example is only a minor modification of the
one given in \cite{BaShCon2000},
but the properties derived here for this example are stronger.

\begin{thm}\label{CounterH}
There exists a (non $\sigma$-compact) subset $H$ of $\Increasing$, such that:
\be
\item $|H|=\b$,
\item $H$ satisfies $\sone(\Gamma,\cO)$ and all finite powers of $H$ satisfy $\ufin(\Gamma,\Gamma)$, but
$H\sm Q$ does not satisfy $\ufin(\Gamma,\Gamma)$; and
\item If $\b=\d$, then $H\sm Q$ does not satisfy $\ufin(\Gamma,\cO)$.
\ee
\end{thm}
\begin{proof}
First, observe that an uncountable
$\sigma$-compact set of reals contains a perfect set,
and this is ruled out by the property $\sone(\Gamma,\cO)$ \cite{coc2}.
The rest of the proof is divided into several lemmas.
We take $H=B\cup Q$, where $B$ is the set
described in the following lemma.

\begin{lem}\label{Blemma}
There exists an unbounded subset $B=\{g_\alpha : \alpha < \b\}$
of $\NN$ where: Each $g_\alpha$ is increasing,
$g_\alpha \leq^* g_\beta$ for each $\alpha<\beta$, and
if $\b=\d$ then $B$ is dominating.
\end{lem}
\begin{proof}
Start with $\{f_\alpha : \alpha<\b\}$ unbounded and $\{h_\alpha : \alpha<\d\}$ dominating,
and by induction on $\alpha$,
choose $g\in\NN$ bounding $\{g_\beta : \beta<\alpha\}$, and if possible
choose also $h\in\NN$ bounding $\{h_\beta : \beta<\alpha\}$ (otherwise take $h\equiv 0$).
Let $g_\alpha$ be an increasing function such that $g,h,f_\alpha\le^* g_\alpha$.
\end{proof}

According to a theorem of Hurewicz (see Rec\l{}aw \cite{RECLAW}),
a zero-dimensional set of reals $X$
satisfies $\ufin(\Gamma,\Gamma)$ (respectively, $\ufin(\Gamma,\cO)$) if, and
only if, each continuous image of $X$ in $\NN$ is bounded (respectively, not dominating).
As $H\sm Q=B$ is unbounded, it does not satisfy $\ufin(\Gamma,\Gamma)$.
If $\b=\d$, then $B$ is dominating, thus $H\sm Q$ does not satisfy $\ufin(\Gamma,\cO)$.

We now show that $H$ satisfies $\sone(\Gamma,\cO)$.
A subset $A$ of $\NN$ is \emph{strongly unbounded}
if for each $f\in\NN$, $|\{g\in A : g\le^* f\}|<|A|$.
Observe that $B$ is strongly unbounded.
For a cardinal $\kappa$, a set of reals $X$ is \emph{$\kappa$-concentrated}
on a set $Y$ if for each open set $U\spst Y$, $|A\sm U|<\kappa$.

\begin{lem}
Assume that $A$ is a strongly unbounded subset of $\NN$ and
$\kappa=|A|$. Then:
\be
\item $A$ is $\kappa$-concentrated on $Q$; and
\item For each family $\scrA$ of open covers of $A\cup Q$, if
$\kappa\le\non(\sone(\scrA,\cO))$, then $A\cup Q$ satisfies
$\sone(\scrA,\cO)$.
\ee
\end{lem}
\begin{proof}
(1) Assume that $U\spst Q$ is open.
Then $\Increasing\sm U$ is a closed and therefore compact subset of
$\Increasing$. Since $\Increasing\sm U$ is disjoint from $Q$, it is a compact
(and therefore bounded) subset of $\NN$.
As $A$ is strongly unbounded, $|A\sm U|=|A\cap(\Increasing\sm U)|<\kappa$.

(2) Assume that $\cU_n\in\scrA$, $n\in\N$. Enumerate $Q=\{q_n :
n\in\N\}$, and choose for each $n$ $U_{2n}\in\cU_{2n}$ such that
$q_n\in U_{2n}$. Let $U = \Union_n U_{2n}$. Then $|A\sm
U|<\kappa\le\non(\sone(\scrA,\cO))$, thus $A\sm U$ satisfies
$\sone(\scrA,\cO)$. For each $n$ choose an element
$U_{2n+1}\in\cU_{2n+1}$ such that $A\sm U\sbst\Union_n U_{2n+1}$.
Then $\seq{U_n}$ is a cover of $A\cup Q$.
\end{proof}

We need the following extension of Lemma 2 of \cite{BaShCon2000}.
\begin{lem}\label{powerlemma}
Assume that $Q^k\sbst X^k\sbst\powInc{k}$, and $\Psi:\to\NN$ is continuous on $Q^k$.
Then there exists $g\in\NN$ such that for each $n$ and each
$x_1,\dots,x_k\in X$,
$$\mbox{if }g(n)<\min\{x_1(n),\dots,x_k(n)\}\mbox{, then }
\Psi(x_1,\dots,x_k)(n)\le g(n).$$
\end{lem}
\begin{proof}
For each $A\sbst\Inc$, let $A\rest n=\{x\rest n : x\in A\}$. For
each $n$, let $\nInc=\Inc\rest n$. For $\sigma\in\nInc$, write
$q_\sigma$ for $q_{\sigma\rest  m}$ where $m=1+\max\{i<n :
\sigma(i)<\infty\}$.

If $\sigma\in\nInc$ and $I$ is a basic open neighborhood of
$q_\sigma$, then there exists a natural number $N$ such that for
each $x\in\Inc$ with $x\rest n\in I\rest n$ and $x(n)>N$, $x\in
I$.

Fix $n$. Use the continuity of $\Psi$ on $Q^k$ to choose, for each
$\vec\sigma = (\sigma_1,\dots,\sigma_k) \in \big(\nInc\big)^k$, a
basic open neighborhood
$$I_{\vec\sigma}= I_{\sigma_1}\x \dots\x I_{\sigma_k}\sbst\powInc{k}$$
of $q_{\vec\sigma}=(q_{\sigma_1},\dots,q_{\sigma_k})$ such that
for all $(x_1,\dots,x_k) \in I_{\vec\sigma}\cap X^k$,
$\Psi(x_1,\dots,\allowbreak
x_k)\allowbreak(n)=\Psi(q_{\vec\sigma})(n)$. For each
$i=1,\dots,k$, choose $N_i$ such that for all $x\in\Inc$ with
$x\rest n\in I_{\sigma_i}\rest n$ and $x(n)>N_i$, $x\in
I_{\sigma_i}$. Define $N(\vec\sigma)=\max\{N_1,\dots,N_k\}$.

The set $I_{\vec\sigma}^{(n)}=\{(x_1\rest  n,\dots, x_k\rest  n) :
(x_1,\dots,x_k)\in I_{\vec\sigma}\}$ is open in
$\big(\nInc\big)^k$ and the family $\{I_{\vec\sigma}^{(n)}:
\vec\sigma\in \big(\nInc\big)^k\}$ is a cover of the compact space
$\big(\nInc\big)^k$. Take a finite subcover
$\{I_{\vec\sigma_1}^{(n)},\dots,I_{\vec\sigma_m}^{(n)}\}$ of
$\big(\nInc\big)^k$. Let $N=\max\{N(\vec \sigma_1),\dots,N(\vec
\sigma_m)\}$, and define
$$g(n)=\max\{N,\Psi(q_{\vec \sigma_1})(n), \dots, \Psi(q_{\vec \sigma_m})(n)\}.$$
For all $x_1,\dots,x_k\in X$, let $i$ be such that $(x_1\rest
n,\dots,x_k\rest n)\in I_{\vec \sigma_i}^{(n)}$. If $x_1(n),\dots,
x_k(n)>N$, then
$\Psi(x_1,\dots,x_k)(n)=\Psi(q_{\vec\sigma_i})(n)\le g(n)$.
\end{proof}

It remains to show that
all finite powers of $H$ satisfy $\ufin(\Gamma,\Gamma)$.
We will show, by induction on $k$, that all continuous images in $\NN$ of
the finite powers $H^k$ of $H$ are bounded.

Assume that $\Psi:H^{k+1}\to\NN$ is continuous.
We may assume that all elements in the image of $\Psi$ are increasing.
Let $g\in\NN$ be (increasing, and) as in Lemma \ref{powerlemma}.
By the Lemma \ref{Blemma},
there exists $\alpha<\b$ such that the set $A=\{n : g(n)<g_{\alpha}(n)\}$
is infinite, and for each $\beta>\alpha$, $A\as\{n : g(n)<g_\beta(n)\}$.
Let $\seq{a_n}$ be an increasing enumeration of $A$, and define
$h\in\NN$ by $h(n) = g(a_n)$.

By Lemma \ref{powerlemma},
for all $\alpha_1,\dots,\alpha_k>\alpha$ and all but finitely many $n$,
$$\Psi(g_{\alpha_1},\dots,g_{\alpha_k})(n)\le \Psi(g_{\alpha_1},\dots,g_{\alpha_k})(a_n)\le g(a_n)=h(n).$$
For each
$f\in\{g_\beta : \beta\le\alpha\}\cup Q$
and each $m=1,\dots,k+1$ define
$\Psi_{m,f}:H^k\to\NN$ by
$$\Psi_{m,f}(x_1,\dots,x_k) = \Psi(x_1,\dots,x_{m-1},f,x_{m+1},\dots,x_k).$$
By the induction hypothesis, the image of each function $\Psi_{m,f}$ is bounded.
As there are less than $\b$ many such functions,
we have that $\Psi[H^{k+1}]$ is bounded.

The proof that $H$ satisfies $\ufin(\Gamma,\Gamma)$ is similar \cite{BaShCon2000}.
This completes the proof of Theorem \ref{CounterH}.
\end{proof}

To state the following corollary, we need some preliminaries.
Consider the collection $\Omega^\gpbl$ of open $\omega$-covers $\cU$ of $X$ such that
there exists a partition $\cP$ of $\cU$ into finite sets such that
for each finite $F\sbst X$ and all but finitely many $\cF\in\cP$,
there exists $U\in\cF$ such that $F\sbst U$.
$\sfin(\Omega,\Omega^\gpbl)$ is strictly stronger than $\ufin(\Gamma,\Gamma)$,
and it is also strictly stronger than $\sfin(\Omega,\Omega)$.
A set $X\sbst\R$ is \emph{meager additive} if for each meager set $M\sbst\R$,
$X+M$ is meager.
If $X$ satisfies $\ufin(\Gamma,\Gamma)$ and
has strong measure zero (both properties follow from $\sone(\Omega,\Omega^\gpbl)$),
then it is meager-additive \cite{NSW}, but $\sone(\Omega,\Omega^\gpbl)$ is
strictly stronger than
being meager additive:
Consider the set $X$ of Theorem \ref{p=c-X}. Then $X$ is meager additive \cite{GM},
and therefore so is $X\sm \fin$, but $X\sm\fin$ does not even satisfy
$\ufin(\Gamma,\cO)$.
Let $\SMZ$ denote the collection of strong measure zero sets of reals.
$\cov(\M)\le\non(\SMZ)$, and strict inequality is consistent.
A set of reals $X$ is \emph{perfectly meager} if for each perfect set $P$, $X\cap P$ is
meager in the relative topology of $P$.
It is \emph{universally meager} if it
does not contain an injective Borel image of a nonmeager set of reals.

\begin{cor}\label{nicecor}
The set $H$ has the following properties:
\be
\item $\sfin(\Omega,\Omega^\gpbl)$ and $\sone(\Gamma,\cO)$,
\item\label{pm}
It is universally meager (in particular, it is perfectly meager); and
\item\label{mgradd}
If $\b\le\non(\SMZ)$, then $H$ satisfies
$\sone(\Omega,\Omega^\gpbl)$. In particular, in this case it is meager-additive.
\ee
\end{cor}
\begin{proof}
(1) $X$ satisfies $\sfin(\Omega,\Omega^\gpbl)$ if, and only if, all finite powers
of $X$ satisfy $\ufin(\Gamma,\Gamma)$ \cite{coc7}.

(2) An uncountable set of reals $X$ satisfying $\sone(\Gamma,\cO)$ cannot
contain a perfect set of reals \cite{coc2}.
Zakrzewski \cite[Proposition 2.3]{zakrUFC} proved
that if $X$ satisfies $\ufin(\Gamma,\Gamma)$ and does not contain
a perfect set, then $X$ is universally meager.

(3) If $\b\le\non(\SMZ)$ then $H$ is $\non(\SMZ)$-concentrated on the countable
set $Q$, which implies that $H$ has strong measure zero. By \cite{prods},
$\sfin(\Omega,\Omega^\gpbl)\cap\SMZ=\sone(\Omega,\Omega^\gpbl)$.
\end{proof}

Corollary \ref{nicecor}(\ref{mgradd}) extends Theorem 2(1)
of \cite{tomekMadd}, which asserts that
if $\b=\aleph_1$ then there exists a meager-additive set of reals.
Corollary \ref{nicecor}(\ref{pm}) implies a negative answer to Steprans'
Question 5 from \cite{stepransQns}:
Does the inequality $\non(\M) > \aleph_1$ imply that no set of size
greater than $\aleph_1$ is perfectly meager? The answer is ``No'', since
$H$ is universally meager, and $\b>\aleph_1$ is consistent with the assumption of
the question.

\begin{prob}\label{Hprob}
~\be
\item Does the set $H$ constructed in Theorem \ref{CounterH} satisfy $\sone(\Gamma,\Gamma)$?
\item Do all finite powers of $H$ satisfy $\sone(\Gamma,\cO)$?
\ee
\end{prob}

The methods of \cite{wqn} may be relevant to Problem \ref{Hprob}(1).
A positive answer to Problem \ref{Hprob}(2) would imply that $H$
satisfies $\sone(\Gamma,\Omega)$.

\medskip

We now treat Menger's Conjecture.
A set of reals $X$ satisfies $\sfin(\Omega,\Omega)$ if, and only if,
all finite powers of $X$ satisfy $\ufin(\Gamma,\cO)$ \cite{coc2}.
We do not know whether there exists in ZFC a (non $\sigma$-compact)
set of size $\d$ satisfying $\sfin(\Omega,\Omega)$.
It is known that $\cov(\M)=\c$ is enough to deduce the existence of
such a set \cite{coc2}.
We will show that this can also be deduced from an assumption which
contradicts $\cov(\M)=\c$.
A subset $F$ of $\NN$ is \emph{finitely-dominating}
if for each $g\in\NN$ there exist $k$ and $f_1,\dots,f_k\in\NN$
such that $g(n)\le^* \max\{f_1(n),\dots,f_k(n)\}$.
Let $\Dfin$ denote the collection of sets
$F$ of increasing elements of $\NN$, which are not finitely-dominating,
and let $\add(\Dfin)=\min\{|\cF| : \cF\sbst\Dfin\mbox{ and }\cup\cF\nin\Dfin\}$.
If $\cov(\M)=\c$ then $\add(\Dfin)=2$, but
NCF is equivalent to $\add(\Dfin)>2$, and
if $\fu<\g$, then $\add(\Dfin)=\c$ \cite{AddQuad}.

\begin{thm}\label{CounterM}
There exists a (non $\sigma$-compact) subset $M$ of $\Increasing$, such that:
\be
\item $|M|=\d$,
\item $M$ satisfies $\sone(\Gamma,\cO)$, but $M\sm Q$ does not satisfy $\ufin(\Gamma,\cO)$,
\item If NCF holds, then $M$ satisfies $\ufin(\Gamma,\Omega)$; and
\item If $\b=\d$
\footnote{Equivalently, a union of less than $\d$ nondominating subsets of $\NN$
is not dominating.}
or $\add(\Dfin)=\d$, then $M$ satisfies $\sfin(\Omega,\Omega)$.
\ee
\end{thm}
\begin{proof}
The proof is similar to that of Theorem \ref{CounterH}. We describe only the differences.
Here, we take $M=D\cup Q$ where $D$ is defined in the following lemma.
\begin{lem}\label{Dlemma}
There exists a dominating subset $D =\{g_\alpha : \alpha < \d\}$
of $\NN$ where
each $g_\alpha$ is increasing, and
for each $f\in\NN$ there exists $\alpha_0<\d$ such that for any
finite set $F\sbst\d\setminus\alpha_0$,
$f(n) < \min\{g_\beta(n) : \beta\in F\}$
for infinitely many $n$.
\end{lem}
\begin{proof}
Let $\{f_\alpha : \alpha<\d\}$ be a dominating subset of $\NN$.
We construct $g_\alpha$ by induction on $\alpha < \d$.
Assume that $g_\beta$ are constructed for $\beta<\alpha$,
such that for each $\beta<\alpha$ and each finite set
$F\sbst\alpha\sm\beta$, the set
$$X_{\beta,F}=\{n : f_\beta(n) < \min\{g_\gamma(n) : \gamma\in F\}\}$$
is infinite.
For each $\beta<\alpha$ and finite $F\sbst\alpha\sm\beta$,
let $h_{\beta,F}\in\NN$ be the increasing enumeration of $X_{\beta,F}$.
The collection
$$\{f_\beta\circ h_{\beta,F} : \beta<\alpha,\ F\text{ a finite subset of }\alpha\sm\beta \}$$
has less than $\d$ many elements, thus there exists an increasing
$g_\alpha\in\NN$ such that $f_\alpha\le^* g_\alpha$, and
for each $\beta$ and $F$,
$g_\alpha\not \le^* f_\beta\circ h_{\beta,F}$, therefore
$g_\alpha\circ h_{\beta,F} \not\le^* f_\beta\circ h_{\beta,F}$,
thus
$$f_\beta(n) < \min\{g_\gamma(n) : \gamma\in F\cup\{\alpha\}\}$$
for infinitely many $n$.
This completes the inductive step.
\end{proof}

We need to prove that
if $\add(\Dfin)>2$, then $M$ satisfies $\ufin(\Gamma,\Omega)$, and
if $\b=\d$ or $\add(\Dfin)=\d$, then $M$ satisfies $\sfin(\Omega,\Omega)$.
The proof of the first assertion uses arguments
similar to the following ones, and we omit it.

Assume that $\add(\Dfin)=\d$ (respectively, $\b=\d$).
Assume that $M^k$ satisfies $\ufin(\Gamma,\Omega)$ (respectively, $\ufin(\Gamma,\cO)$).
Assume that $\Psi:M^{k+1}\to\NN$ is continuous and all elements in its image are increasing.
By \cite{huremen1} (respectively, by Hurewicz' Theorem),
it suffices to show that $\Psi[M^{k+1}]$
is not finitely dominating (respectively, not dominating).

Let $g\in\NN$ be increasing and as in Lemma \ref{powerlemma},
and take $\alpha_0<\d$ as in Lemma \ref{Dlemma}.
Assume that
$F=\{g^i_j : i<m, j\le k\}\sbst\{g_\alpha : \alpha\in\d\setminus\alpha_0\}$.
Then there exists an infinite set $A\sbst\N$ such that
$g(n)<\min\{g^i_j(n) : i<m, j\le k\}$ for each $n\in A$,
and if $y_i=\Psi(g^i_0,\dots,g^i_k)$ for each $i<m$, then
by Lemma \ref{powerlemma},
$$\max\{y_0(n),\dots,y_{m-1}(n)\}\le g(n)$$
for each $n\in A$.
Thus, $\Psi[\{g_\alpha : \alpha\in\d\setminus\alpha_0\}^{k+1}]$ is not
finitely dominating.
It follows that $\Psi[M^{k+1}]$ is a union of less that $\d$ many sets which
are not finitely dominating.
%
\end{proof}

\begin{rem}
If one only wants to obtain (1) and (2) of Theorem \ref{CounterM}, then
it suffices to take $M=\Psi[D]\cup\Q$ where $D$ is any strongly
unbounded subset of $\NN$ and $\Psi:\NN\to\R\sm\Q$ is a homeomorphism.
Essentially, this idea goes back to Rothberger.
\end{rem}

\paragraph{\textbf{Acknowledgements}}
We thank Roman Pol, Taras Banakh, and Lubomyr Zdomsky for making
crucial comments on earlier versions of this paper,
and the referee for extending and simplifying our proof of
Proposition \ref{referee}. We also thank Tomasz Weiss for fruitful
discussions.

Chaber and Pol proposed another approach
to the Menger and Hurewicz conjectures \cite{ChaPol}.

\end{document}